\documentclass[a4paper,8pt]{article}
\begin{document}
\bf

\centerline{Matching on a line}
	
\bigskip
\bf
\centerline{Josef Bukac}

\bigskip
\rm	
\centerline{Bulharska 298}
\centerline{55102 Jaromer-Josefov}
\centerline{Czech Republic}

\bigskip
\bf
\noindent
Keywords:
\rm 
bipartite, design of experiments, nonbipartite, $n$-tuple, triple, 
tripartite, quadripartite.

\bigskip
\bf
\noindent
MSC2010:
\rm
90C27, 68Q17, 05C85.	
	
\bigskip
\bf
\noindent
Abstract:
\rm
Matching is a method of the design of experiments. If we had 
an even number of patients and wanted to form pairs of 
patients such that their ages, for example,  in each pair 
be as close as possible, we would use nonbipartite matching.
Not only do we present a fast method to do this, we also
extend our approach to triples, quadruples, etc.

In part 1 a matching algorithm uses $kn$ points on a line as vertices,
pairs of vertices as edges, and either absolute values of differences 
or the squares of differences as weights
or distances. It forms $n$ of $k$-tuples with the minimal sum of 
distances within each $k$-tuple in $O(n$  $log$ $n)$ time. 

In part 2 we present a trivial algorithm for bipartite matching 
with absolute values or squares of differences as weights and 
a generalisation to tripartite matching on tripartite graphs.

\bigskip
\bf
Introduction

\bigskip
\rm
Further references about the use of nonbipartite matching in experimental design are in survey papers Beck (2015) or Lu (2011). 
Imagine we have something like 300 patients
and we want to form pairs of patients such that the ages, 
as an example of a confounding variable, of the patients in each 
pair are as close as possible. 
Our goal is to make applications of treatment A and treatment B comparable. In our paper we also show how to form triples that 
are convenient for an application of placebo, treatment A, or 
treatment B. Quadruples may be used for an application of 
placebo, treatment A,  treatment B, and interaction of A and B combined.

We picked age as an example of a trivial confounding 
variable but there are sofisticated ways of defining such 
a variable, propensity score being one of them.

Matching is used when we want to avoid the
effect of a confounding variable. If there are more
such variables, it is customary to agregate them 
to get just one variable typically called 
a scale or score. Applications vary in the fields of 
medicine, social sciences, psychology, or education. 

Should we want to use an $n$-dimensional space for 
$n$ confounding variables we would have to multiply 
each of these variables by some constant to take care 
of their importance, units, etc, and we would have 
to derive those constants somehow only to find out 
that scores are a better choice.

The repeatability of matching is important because, 
unlike randomization, matching gives the same result 
each time it is repeated, save for ties.

In such a setting, each individual becomes a vertex 
of a simple complete graph and the weight of each edge 
is defined as the absolute value $A$ of the difference 
between their scores or it is defined as the square $S$
of this difference.

In part 1 we want to show that the calculation of 
nonbipartite matching becomes trivial. Also triangle 
matching, termed $3$-matching, becomes an easy task and 
so does $4$-matching, and generally $n$-matching so far 
for $n\le 16$ in the case of absolute values of differences 
as weights or $n\le 8$ when the sum of squares of differences 
is used .

We consider a complete simple graph $G$ with an even 
number $|V|$ of vertices $V.$ The set of edges is denoted 
as $E.$ Let $M$ be a subset of $E.$ $M$ is called 
a matching if no edges in $M$ are adjacent in $G.$  
A matching $M$ is called perfect if each vertex in $V$
is incident to some edge in $M.$ We assume a weight $w(e)\ge 0$
is attached to every edge $e\in E.$ We are then looking
for a perfect matching $M$ for which the sum of weights 
is maximal. Equivalently, we may look for a perfect matching
with a minimal sum of weights by picking some upper bound $u$
of weights and form new weights as $u-w(e)$ for each $e\in E.$ 

The Edmonds method for finding a maximal weighted matching 
in a general weighted graph is presented in Papadimitriou (1982). 
The time necessary for calculation is polynomial, $O(|V|^3)$, in the 
number of vertices. Writing the program would be tedious 
but we recommend Beck (2015) or an internet address through which 
the solution may be obtained:

\noindent
http://biostat.mc.vanderbilt.edu/wiki/Main/NonbipartiteMatching

Even though the running time is polynomial, the degree 3 may
turn out to be too high for practical calculations for a large number
of vertices $|V|$ requiring a large $|V|$ by $|V|$ matrix of distances.

If the number of vertices is divisible by three and a constant 
$B$ is given, the decision problem whether there is a 
perfect 3-matching such that the sum over all the
triples of the distances between the three points in each 
triple is less than some constant $B$ is known to be $NP$-complete.
We may refer to the problem exact cover by 3-sets in Garey (1979)
or Papadimitriou (1982). That is why the problems we study seem 
to be so discouraging.

In the second part of the paper we show a trivial method of calculating 
a minimal perfect matching on a regular complete bipartite graph with weights 
on edges being absolute values of differences of weights on vertices.
This is only a stepping stone to the design of a method of calcualtion 
of a minimal perfect matching on a regular complete tripartite graph
with the same definition of weights.

\bigskip
\bf
Part 1

\bigskip
\bf
1.1. Matching on a line

\bigskip
\rm
We study a complete graph the vertices of which are
points on a real $x$-axis. 
These points are denoted as $x_i$.

The edges are the line segments between these points. The weight
of each edge $(x_i, x_j)$ is defined as the distance 
of its two endpoint $|x_i-x_j|.$ 

Since the possibility of repeated observations is common 
in statistics, we do not use sets, we use the notion of 
a $k$-tuple. We call $(1, 1, 2)$ a triple whereas a set 
would consist of two elements, 1 and 2.

\bigskip
\bf
\noindent
Definition 1.1.
\rm\quad
We define a distance $A$ within a $k$-tuple as the sum of the 
distances of all the pairs formed of the elements of
the $k$-tuple. 
$$A(x_1, x_2,\dots , x_k)=
\sum_{i=1}^{k-1}\sum_{j=i+1}^{k}|x_j-x_i|.$$

\bigskip
\rm
There are $k(k-1)/2$ summands in this formula.
Our calculations will be simplified by the following.

\bigskip
\bf
\noindent
Definition 1.2.
\rm
A $k$-tuple $(x_1, x_2,\dots ,x_k)$ is sorted
if $x_1\le x_2\le\dots\le x_k.$

\bigskip
\bf
\noindent
Theorem 1.1.
\rm
If the $k$-tuple is sorted, the distance within the $k$-tuple
is 
$$A(x_1, x_2,\dots , x_k)=\sum_{i=1}^k (2i-k-1)x_i$$

\bf
\noindent
Proof.\rm\quad 
$$A(x_1, x_2,\dots , x_k)
=\sum_{i=1}^{k-1}\sum_{j=i+1}^{k}(x_j-x_i)
=\sum_{i=1}^{k-1}\sum_{j=i+1}^{k}x_j-\sum_{i=1}^{k-1}\sum_{j=i+1}^{k}x_i$$ 
$$=x_2+2x_3+\dots +(k-1)x_k-\sum_{i=1}^{k-1}(k-i)x_i
=\sum_{i=1}^k (i-1)x_i-\sum_{i=1}^{k-1}(k-i)x_i$$

\noindent
$\quad =\sum_{i=1}^k (2i-k-1)x_i.$

\bigskip
\rm
For example, if a sorted pair $(x_1, x_2), $ $x_1\le x_2,$ is given, 
the distance is defined as $A(x_1, x_2)=x_2-x_1.$ 
For a sorted triple $(x_1, x_2, x_3),$ $x_1\le x_2\le x_3,$ we define 
the distance within as $A(x_1, x_2, x_3)=2(x_3-x_1).$

\bigskip
\bf
\noindent
Definition 1.3.
\rm
Let a $kn$-tuple be given. A partition of this $kn$-tuple
into $n$ of $k$-tuples is called a $k$-tuple partition
of a $kn$-tuple.

\bigskip
\bf
\noindent
Definition 1.4.
\rm
Let a $kn$-tuple be given. A $k$-tuple partition of 
this $kn$-tuple is called minimal if the sum of
the distances within taken over all the $n$ of $k$-tuples is 
less than or equal to the sum of distances within taken
over $k$-tuples of any other $k$-tuple partition.

\bigskip
\rm
We may assume there may be more than one minimal partition.
For a $kn$-tuple there are $(kn)!\big/(k!)^n$ $k$-tuple
partitions from which we want to find those with minimal sum of 
distances within $k$-tuples. They are too many for 
the brute force method to work for a large $n.$ But
it will work for $n=2$ if $k$ is small.

The greedy method will not work either. That can be shown 
by way of an example $(1, 3, 4, 5, 8, 9)$ in which the
triple with the smallest sum of distances within is $(3, 4, 5),$ 
$A(3,4,5)=4$, the distance within the remaining items
is $A(1,8,9)=16.$ The sum is $A(3,4,5)+A(1,8,9)=20.$
We get a smaller sum of distances within  if we take $(1, 3, 4)$ and 
$(5, 8, 9)$ yielding $A(1,3,4)+$ $A(5,8,9)=$ 
$6+8=$ $14.$

First we want to show what a minimal partition 
for $2k$-tuples looks like. If we can do that, we will 
use induction to show it works for any $n>2.$

\bigskip
\bf
\noindent
Theorem 1.2.
\rm
Let a sorted $4$-tuple $(x_1, x_2, x_3, x_4)$ be given.
Then the two pairs $(x_1, x_2)$ and $(x_3, x_4)$ 
have a minimal sum of distances defined as absolute 
values of differences.

\bigskip
\bf
\noindent
Proof.\rm\quad
The sum of distances of $(x_1, x_2)$ and $(x_3, x_4)$ is
$x_2-x_1+x_4-x_3.$
We form other possible partitions into pairs, calculate 
the sum of their distances, and compare them with the sum 
of distances of $(x_1, x_2)$ and $(x_3, x_4).$ 

Other possible sorted pairs are: 

1) $(x_1, x_3)$ and $(x_2, x_4).$ The difference is
$$A(x_1, x_3)+A(x_2, x_4)-(A(x_1, x_2)+A(x_3, x_4))
=2(x_3-x_2)\ge 0.$$

2) $(x_1, x_4)$ and $(x_2, x_3).$ The difference is
$$A(x_1, x_4)+A(x_2, x_3)-(A(x_1, x_2)+A(x_3, x_4))
=2(x_3-x_2)\ge 0.$$

\bigskip
\bigskip
\bf
\noindent
Theorem 1.3.
\rm
Let a sorted $6$-tuple $(x_1, x_2,\dots x_6)$ be given.
Then the two sorted triples $(x_1, x_2, x_3)$ and 
$(x_4, x_5, x_6)$ have a minimal sum of distances within
defined as the sum of absolute values of differences.

\bigskip
\bf
\noindent
Proof.\rm\quad
We form all the other possible triples, calculate the sum of their
distances within, and subtract from them the sum of distances
$A(x_1, x_2, x_3 )$ $+A(x_4, x_5, x_6).$ 

Other possible triples are listed in such a way that $x_1$
appears in the first triple because the order does not matter
in this case: 

\noindent
1)\quad $(x_1, x_2, x_4)$ and $(x_3, x_5, x_6)$

\noindent
$A(x_1, x_2, x_4)+A(x_3, x_5, x_6)-A(x_1, x_2, x_3)-A(x_4, x_5, x_6)=$

$4(x_4-x_3)\ge 0$

\smallskip
\noindent
2)\quad $(x_1, x_2, x_5)$ and $(x_3, x_4, x_6)$

\noindent
$A(x_1, x_2, x_5)+A(x_3, x_4, x_6)-A(x_1, x_2, x_3)-A(x_4, x_5, x_6)=$

$2(x_5+x_4-2x_3)\ge 0$

\smallskip
\noindent
3)\quad $(x_1, x_2, x_6)$ and $(x_3, x_4, x_5)$

\noindent
$A(x_1, x_2, x_6)+A(x_3, x_4, x_5)-A(x_1, x_2, x_3)-A(x_4, x_5, x_6)=$

$2(x_5+x_4-2x_3)\ge 0$

\smallskip
\noindent
4)\quad $(x_1, x_3, x_4)$ and $(x_2, x_5, x_6)$

\noindent
$A(x_1, x_3, x_4)+A(x_2, x_5, x_6)-A(x_1, x_2, x_3)-A(x_4, x_5, x_6)=$

$2(2x_4-x_3-x_2)\ge 0$

\smallskip
\noindent
5)\quad $(x_1, x_3, x_5)$ and $(x_2, x_4, x_6)$

\noindent
$A(x_1, x_3, x_5)+A(x_2, x_4, x_6)-A(x_1, x_2, x_3)-A(x_4, x_5, x_6)=$

$2(x_5+x_4-x_3-x_2)\ge 0$

\smallskip
\noindent
6)\quad $(x_1, x_3, x_6)$ and $(x_2, x_4, x_5)$

\noindent
$A(x_1, x_3, x_6)+A(x_2, x_4, x_5)-A(x_1, x_2, x_3)-A(x_4, x_5, x_6)=$

$2(x_5+x_4-x_3-x_2)\ge 0$

\smallskip
\noindent
7)\quad $(x_1, x_4, x_5)$ and $(x_2, x_3, x_6)$

\noindent
$A(x_1, x_4, x_5)+A(x_2, x_3, x_6)-A(x_1, x_2, x_3)-A(x_4, x_5, x_6)=$

$2(x_5+x_4-x_3-x_2)\ge 0$

\smallskip
\noindent
8)\quad $(x_1, x_4, x_6)$ and $(x_2, x_3, x_5)$

\noindent
$A(x_1, x_4, x_6)+A(x_2, x_3, x_5)-A(x_1, x_2, x_3)-A(x_4, x_5, x_6)=$

$2(x_5+x_4-x_3-x_2)\ge 0$

\smallskip
\noindent
9)\quad $(x_1, x_5, x_6)$ and $(x_2, x_3, x_4)$

\noindent
$A(x_1, x_5, x_6)+A(x_2, x_3, x_4)-A(x_1, x_2, x_3)-A(x_4, x_5, x_6)=$

$2(2x_4-x_3-x_2)\ge 0$

That finishes the proof. It was actually generated by a comuter.

\bigskip
\rm
We may consider a pair of $k$-tuples for any $k\ge 2$ 
but we have to generate all such pairs of $k$-tuples 
while keeping $x_1$ in the first of them. It means 
the number of all such pairs is ${2k-1}\choose {k-1}$. 
We do not need any symbolic algebra to do this 
for it will suffice to keep in mind that the $k$-tuples
are generated as combinations represented as subscripts. 
The sums of the two distances within each $k$-tuple are 
expressed as coefficients assigned to subscripts. 
The resulting inequality is obtained by comparing 
the coefficients as indicated in Theorems 1.2 and 1.3. 

Not only does the amount of work on each $k$-tuple increase 
approximately proportionally with respect to $k,$
but, and more importantly, the growth of binomial coefficients 
${2k-1}\choose {k-1}$ becomes prohibitive for  
calculations as $k$ increases. When we have the
time of calculations for $k,$ the time necessary
for $k+1$ will be approximately equal to the time for $k$ 
times the following

$$\frac{k+1}{k}{{2(k+1)-1}\choose{(k+1)-1)}}\big/
{{2k-1}\choose{k-1)}}=4+\frac{2}{k}$$
This is the reason why we have been able to verify our claims 
so far for $k\le 16.$ We stopped at 16 also because it can be 
used to form $4$ by $4$ tables.

\bigskip
\bf
\noindent
Definition 1.5.
\rm
Let $(x_1,x_2,\dots ,x_k)$ and $(y_1,y_2,\dots ,y_k)$ be
two distinct sorted $k$-tuples. The smallest subscript $j$ 
for which $x_j\ne y_j$ is called the smallest subscript
of discordance.

\bigskip
We note that if $(x_1,x_2,\dots ,x_k)=(y_1,y_2,\dots ,y_k),$
the smallest subscript of discordance is not defined.

\bigskip
\bf
\noindent
Theorem 1.4.
\rm
If for any sorted $2k$-tuple  $(x_1, x_2,\dots , x_{2k})$ 
the two sorted $k$-tuples $(x_1, x_2,\dots , x_{k})$ and 
$(x_{k+1}, x_{k+2},\dots , x_{2k})$ are the minimal solution
of the $k$-matching problem, then the minimal solution
of the $k$-matching problem for a sorted $kn$-tuple, $n>0,$ is 
given by $n$ sorted $k$-tuples 
$$(x_{(i-1)k+1}, x_{(i-1)k+2},\dots , x_{(i-1)k+k})$$ 
for $i=1,\dots , n.$

\bigskip
\bf
\noindent
Proof.
\rm\quad
We prove the theorem by induction. If $n=1,$ the theorem is
obvious. If $n=2,$ the theorem follows directly from its
assumption. We assume the theorem is true if $n-1>1$ and show
it is true for $n.$

We consider all the possible minimal $k$-tuple partitions.
If there is an $k$-tuple partition such that for some sorted 
$k$-tuple $(y_1,y_2,\dots ,y_k)$ we have $x_i=y_i$ for all 
$i=1, 2,\dots ,k,$ we are done and we may also exclude
the case that the smallest subscript of discordance is not defined
in the following.

If $(x_1,x_2,\dots ,x_k)\ne (y_1,y_2,\dots ,y_k),$ we will show a contradiction.
We will compare $k$-tuples with $(x_1, x_2, \dots , x_k).$
Out of all the minimal $k$-tuple partitions we pick the one 
containing the $k$-tuple $(y_1, y_2, \dots , y_k)$ 
for which the smallest subscript of discordance $j$ 
is the highest. It is obvious that for such a $k$-tuple 
$y_1=x_1$ holds for otherwise the smallest subscript of 
discordance $j$ would be $1.$ If $y_1=x_1,$ we have $j>1.$ 

Since $j,$ where $1<j\le k,$ is the lowest subscript for which
$y_j\ne x_j,$  this $x_j$ must be in some other $k$-tuple
$(z_1, z_2, \dots , z_k)$ in the partition. We concatanate 
these two $k$-tuples to obtain a $2k$-tuple 
$(y_1, y_2, \dots , y_k, z_1, z_2, \dots , z_k)$
and apply the assumption of the theorem to obtain a minimal 
solution of the $k$-matching problem on this $2k$-tuple as 
$(t_1, t_2,\dots , t_j,\dots , t_k )$
and  $(u_1, u_2,\dots , u_k )$ where
$t_i=x_i$ for $i=1, 2, \dots , j.$

We have two cases. We obtain a $k$-tuple 
partition containing $(x_1,x_2, \dots ,x_k)$ which is 
a contradiction. 

If we do not obtain a partition containing 
$(x_1,x_2, \dots ,x_k),$  the smallest
subscript of discordance is $i,$ where $j<i<k,$
when  
$(t_1, t_2,\dots , t_k)$ is compared with 
$(x_1, x_2,\dots , x_k).$
We have obtained a $k$-tuple partition for which 
the smallest subscript of discordance is higher
than $j,$ contradicting the assumption. 

The proof is finished by removing $x_1,x_2,\dots ,x_k$
from the original $kn$-tuple obtaining a $(n-1)k$-tuple.

\bigskip
\bf
\noindent
Corollary.
\rm
If the assumption in theorem 1.4 holds,  then the minimal solution
of the $k$-matching problem for a not necessarily sorted 
$kn$-tuple, $n>0,$ is obtained in the running time necessary 
for sorting the $kn$-tuple.

\bigskip
\rm
It means the matching problem is solved in $O(N\log N)$ time where
$N=kn$ is the number of items to be matched. 

\bigskip
\bf
1.2 Sum of squares of differences

\bigskip
\rm
We all know that statisticians would prefer the sum of squares
of all differences to evaluate the distance within a $k$-tuple.
Let $(x_1, x_2,\dots , x_k)$ be a $k$-tuple the distance within
will be defined as
$$S(x_1, x_2,\dots , x_k)=
\sum_{i=1}^{N-1}\sum_{j=i+1}^N(x_j-x_i)^2$$
Some avid statisticians would even require the minimization of
the sum of variances but this is equivalent to the sum of squares 
of all the differences as explained in the appendix.

First we want to show what a minimal partition for $2k$-tuples 
looks like.

\bigskip
\bf
\noindent
Theorem 1.5.
\rm
Let a sorted $4$-tuple $(x_1, x_2, x_3, x_4)$ be given.
Then the two pairs $(x_1, x_2)$ and $(x_3, x_4)$ 
have a minimal sum of squares of all differences.

\bigskip
\bf
\noindent
Proof.\rm\quad
The sum of distances of $(x_1, x_2)$ and $(x_3, x_4)$ is
$S(x_1, x_2)+S(x_3, x_4)=(x_2-x_1)^2+(x_4-x_3)^2.$
We form other possible partitions into pairs, calculate 
the sum of the differences squared, and compare them with the sum 
of distances $S(x_1, x_2)+S(x_3, x_4).$ 
Other possible sorted pairs are: 

\smallskip
\noindent
1) $(x_1, x_3)$ and $(x_2, x_4).$ The difference is

$S(x_1, x_3)+S(x_2, x_4)-(S(x_1, x_2)+S(x_3, x_4))
=2(x_4-x_1)(x_3-x_2)\ge 0.$

\smallskip
\noindent
2) $(x_1, x_4)$ and $(x_2, x_3).$ The difference is

$S(x_1, x_4)+S(x_2, x_3)-(S(x_1, x_2)+S(x_3, x_4))
=2(x_3-x_1)(x_4-x_2)\ge 0.$

\bigskip
\bigskip
\bf
\noindent
Theorem 1.6.
\rm
Let a sorted $6$-tuple $(x_1, x_2,\dots x_6)$ be given.
Then the two sorted triples $(x_1, x_2, x_3)$ 
and $(x_4, x_5, x_6)$ have a minimal sum of squares
of all differences.

\bigskip
\bf
\noindent
Proof.\rm\quad
We form all the other possible triples, calculate the sum of their
distances, and subtract from them the sum of distances
$S(x_1, x_2, x_3 )+S(x_4, x_5, x_6).$ 

Other possible triples are listed in such a way that $x_1$
appears in the first triple because the order does not matter in this case:

\smallskip
\noindent
1)$S(x_1, x_2, x_4)+S(x_3, x_5, x_6)-S(x_1, x_2, x_3)-S(x_4, x_5, x_6)=$

$2(x_4-x_3)(x_6+x_5-x_2-x_1)\ge 0$

\smallskip
\noindent
2)$S(x_1, x_2, x_5)+S(x_3, x_4, x_6)-S(x_1, x_2, x_3)-S(x_4, x_5, x_6)=$

$2(x_3-x_5)(x_1+x_2-x_4-x_6)\ge 0$

\smallskip
\noindent
3)$S(x_1, x_2, x_6)+S(x_3, x_4, x_5)-S(x_1, x_2, x_3)-S(x_4, x_5, x_6)=$

$2(x_6-x_3)(x_5+x_4-x_2-x_1)\ge 0$

\smallskip
\noindent
4)$S(x_1, x_3, x_4)+S(x_2, x_5, x_6)-S(x_1, x_2, x_3)-S(x_4, x_5, x_6)=$

$2(x_4-x_2)(x_6+x_5-x_3-x_1)\ge 0$

\smallskip
\noindent
5)$S(x_1, x_3, x_5)+S(x_2, x_4, x_6)-S(x_1, x_2, x_3)-S(x_4, x_5, x_6)=$

$2(x_5-x_2)(x_6+x_4-x_3-x_1)\ge 0$

\smallskip
\noindent
6)$S(x_1, x_3, x_6)+S(x_2, x_4, x_5)-S(x_1, x_2, x_3)-S(x_4, x_5, x_6)=$

$2(x_6-x_2)(x_5+x_4-x_3-x_1)\ge 0$

\smallskip
\noindent
7)$S(x_1, x_4, x_5)+S(x_2, x_3, x_6)-S(x_1, x_2, x_3)-S(x_4, x_5, x_6)=$

$2(x_6-x_1)(x_5+x_4-x_3-x_2)\ge 0$

\smallskip
\noindent
8)$S(x_1, x_4, x_6)+S(x_2, x_3, x_5)-S(x_1, x_2, x_3)-S(x_4, x_5, x_6)=$

$2(x_5-x_1)(x_6+x_4-x_3-x_2)\ge 0$

\smallskip
\noindent
9)$S(x_1, x_5, x_6)+S(x_2, x_3, x_4)-S(x_1, x_2, x_3)-S(x_4, x_5, x_6)=$

$2(x_4-x_1)(x_6+x_5-x_3-x_2)\ge 0$

That finishes the proof. 

It is interesting to see that the factorization of all the quadratic forms could
be done. We actually wrote a program that checked the faktorization
for $2\le k\le 8$ and verified the nonnegativity of each factor.

The final step is the use of theorem 1.4 to show how to calculate $k$-matching
for $2\le k\le 8.$ Now we see that it does not matter which of the two
mentioned distances within, either the sum of absolute values of all the differences
or the sum of squares of all the differences, we use, we get the same 
$k$-matching.

\bigskip
\bigskip
\bf
1.3 Statistical applications

\bigskip
\rm
The $k$-tuples obtained by our algorithm are sorted. 
That could have an unpleasant effect on statistical
procedures because of the inequality of the means of the first 
entries of the $k$-tuples as compared with the means of the last 
entries of the $k$-tuples. We know they are different, as long as 
$x_1,x_2,\dots ,x_{kn}$ are not all equal.

We want to avoid randomization in the spirit of our paper. 
What we are trying to achieve is the rearrangement of the items in 
$k$-tuples in such a way that all the means over the $i$-th 
items, $i=1,\dots, k,$ are as close as possible. The minimization 
process reminds us of an NP-complete optimization partition problem 
even though typically the numbers $x_1,x_2,\dots ,x_{kn}$ are not
integers.

Even though partitioning is beyond the scope of this paper, one 
heuristic way of handling the problem is sorting the $k$-tuples
with respect to the distances within in a discending order and 
keeping track of the subtotals starting from the first $k$-tuple 
and rearranging each consecutive $k$-tuple to keep the differences 
as small as possible at each step. This approach obviously does 
not guarantee we obtain the smallest possible differences among 
the means.

\bigskip
\bf
Part 2

\bigskip
\bf
2.1 Bipartite graphs

\bigskip
\rm
Even though algoritms for finding optimal bipartite matching
are so well known that they are presented in introductory textbooks,
such as Bondy(1976), we present another approach because it will 
find applications in tripartite matching.

The regular complete bipartite graph consists of two disjoint
vertex sets $A$ and $B,$ $|A|=|B|=n,$ 
and edges $A\times B.$  We assume that to each of the 
vertices in $A$ and in $B$ respectively a real numbers 
$x_i$ and $y_i$ are assigned as their values. 
The weight associated with each edge $(a_i, b_j)$ is defined as 
either $w_{abs}(a_i, b_j)=|x_i-y_j|$  for all $1\le i, j\le n$ or 
$w_{sq}(a_i, b_j)=(x_i-y_j)^2$  for all $1\le i, j\le n.$

\bigskip
\bf
\noindent
Definition 2.1.
\rm
A perfect matching in a regular complete bipartite graph with vertex set $A\cup B$
is a subset of edges $M_{A,B}$ such that each vertex of $A$ is connected by an 
edge to one vertex of $B$ and each vertex of $B$ is connected 
to one vertex of $A.$

\bigskip

We use the notation $M_{A,B}$ to indicate that we are dealing with the vertex
set  $A\cup B.$

\bigskip
\bf
\noindent
Definition 2.2.
\rm
The weight of a perfect matching is 
$$w(M_{A,B})=\sum_{(a,b)\in M_{A,B}} w(a,b)$$

\bigskip
\bf
\noindent
Definition 2.3.
\rm
A perfect matching $M_{A,B}^{min}$ is minimal if its weight is less than 
or equal to that of any other perfect matching, 
$w(M_{A,B}^{min})\le w(M_{A,B}).$

\bigskip
To avoid any trouble, we mention we do not make any distinction
between the vertices and the values they are assigned, 
we again consider $n$-tuples of real numbers. Sorted $n$-tuples
are described in definition 1.2. 

\bigskip
\bf
\noindent
Theorem 2.1.
\rm
Let two sorted $n$-tuples, $(x_1, x_2, \dots x_n)$ and 
$(y_1, y_2, \dots y_n)$ be given. If the weight of each
edge is defined as the absolute value of the difference
between $x_i$ and $y_j,$ $w_{abs}(a_i,b_j)=|x_i-y_j|,$
for each $1\le i, j\le n,$ then the minimal
perfect matching consists of edges with values 
$(x_1, y_1),$ $(x_2, y_2),$ $\dots ,$ $(x_n, y_n).$

\bigskip
\bf
\noindent
Proof.
\rm\quad
The theorem is true for $n=1.$ If $n>1,$ we assume it is true for $n-1.$ 
In a minimal perfect matching there is an edge with one endpoint value $x_1.$
If the other endpoint value of this edge is $y_1,$ we are done. If not, the other
endpoint value is $y_j$ for some $j>1.$ Another edge must have $y_1$ as its endpoint value,
this edge has some $x_i,$ $i>1$ as its other endpoint value.

Let the required $x_1, x_i, y_1, y_j$ be given
We first consider the three cases when $x_1$ is less
than the rest of the points, $x_i,$ $y_1,$ $y_j.$ Three cases are listed depending on the
position of $x_i.$ 

Two options are possible in each case. The one containing $x_1, y_1$
is subtracted from the other one.  

1. Let  $x_1\le x_i\le y_1\le y_j.$

We subtract $|x_1-y_1|+|x_i-y_j|=$ $y_1-x_1+y_j-x_i$ from
$|x_1-y_j|+|x_i-y_1|=$ $y_j-x_1+y_1-x_i.$ The result is zero
because we subtract the same expression.

2. Let $x_1\le y_1\le x_i\le y_j.$ Then
$|x_1-y_1|+|x_i-y_j|=$ $y_1-x_1+y_j-x_i$ is subtracted from
$|x_1-y_j|+|x_i-y_1|=$ $y_j-x_1+x_i-y_1,$ the difference
is $y_j-x_1+x_i-y_1-(y_1-x_1+y_j-x_i)=$
$2x_i-2y_1=$ $2(x_i-y_1)\ge 0.$

3. Let $x_1\le y_1\le y_j\le x_i.$
Then $|x_1-y_1|+|x_i-y_j|=$ $y_1-x_1+x_i-y_1$ is subtracted from
$|x_1-y_j|+|x_i-y_1|=$ $y_j-x_1+x_i-y_1,$ the difference 
is $y_j-x_1+x_i-y_1-$ $(y_1-x_1+x_i-y_1)=$ 
$y_j-x_1+x_i-y_1-y_1+x_1-x_i+y_1=$
$y_j-y_1\ge 0.$

In the case that $y_1$ is the smallest number we just
swap $x$'s with $y$'s.

We conclude that the minimal matching contains an edge with
endpoint values $(x_1, y_1)$ 
and induction makes sense.

\bigskip
\bf
\noindent
Theorem 2.2.
\rm
Let two sorted $n$-tuples, $(x_1, x_2, \dots x_n)$ and 
$(y_1, y_2, \dots y_n)$ be given. If the weight of each edge
is defined as $w_{sq}=(x_i-y_j)^2,$ for each $1\ge i, j\le n,$ then the minimal
perfect matching consists of edges with values 
$(x_1, y_1),$ $(x_2, y_2),$ $\dots ,$ $(x_n, y_n).$

\bigskip
\bf
\noindent
Proof.
\rm\quad
The theorem is true for $n=1.$ If $n>1,$ we assume it is true for $n-1.$ 
In a minimal perfect matching there is an edge with one endpoint value $x_1.$
If the other endpoint value of this edge is $y_1,$ we are done. If not, the other
endpoint value is $y_j$ for some $j>1.$ Another edge must have $y_1$ as its endpoint value,
this edge has some $x_i,$ $i>1$ as its other endpoint value.

We subtract $(x_1-y_1)^2+(x_i-y_j)^2=x_1^2+y_1^2-2x_1y_1+x_i^2+y_j^2-2x_iy_j$
from $(x_1-y_j)^2+(x_i-y_1)^2=x_1^2+y_j^2-2x_1y_j+y_1^2+x_i^2-2x_iy_1.$
The difference is $-2x_1y_j-2x_iy_1+2x_1y_1+2x_iy_j=2(x_i-x_1)(y_j-y_1)\ge 0.$

\bigskip
\smallskip
The property described in theorems 2.1 or 2.2 not only allows us to
calculate minimal matching quickly, it will be used in the construction
of tripartite matching. The following definition will allow us to formulate
the results in a bit more general but simple setting.

\bigskip
\bf
\noindent
Definition 2.4.
\rm\quad 
Let two sorted $n$-tuples, $(x_1, x_2, \dots x_n)$ and 
$(y_1, y_2, \dots y_n)$ be given. Let a weight of each edge
be defined as $w(x_i,y_j),$ for each $1\ge i, j\le n,$ if the minimal
perfect matching consists of edges with values 
$(x_1, y_1),$ $(x_2, y_2),$ $\dots ,$ $(x_n, y_n)$ for any 
$(x_1, x_2, \dots x_n)$ and $(y_1, y_2, \dots y_n),$ then the weight $w$
is called line matching or LM.

\bigskip
Counterexample: Let a weight be defined as a product $w_p(x_i, y_j)=x_iy_j.$
If we consider $x=(1, 2, 3)$ and $y=(1, 2, 3),$ as an example, 
then the sum of products is $1*1+2*2+3*3=14.$ 
If we use the reverse order $z=(3,2,1),$ then the sum of products is $1*3+2*2+3*1=10<14.$ As a result we can say
that the weight $w_p$ defined as a product is not LM.

We will not study which weights are LM and which are not. 
It suffices to see that the weights $w_{abs}$ and $w_{sq}$ 
are the ones with LM property and those are the ones that would 
be used in practice.

\bigskip
\bf
2.2 Tripartite graphs

\bigskip
\rm
A regular complete tripartite graph is the union of three disjoint
vertex sets $A,$ $B,$ and $C,$ for which $|A|=|B|=|C|=n,$ 
and edges in $A\times B,$ $B\times C,$ and $C\times A.$ 
For any $a\in A,$ $b\in B,$ and  $c\in C$ the edges are denoted as 
$(a,b),$ $(b,c),$ and $(c,a)$ respectively.

\bigskip
We define a matching $M_{A,B,C}\subset A\times B\times C$ as a set of triples of vertices in $M_{A,B,C}$ 
such that if for any two distinct triples $(a_{i_1}, b_{j_1}, b_{k_1})\in M_{A,B,C}$ 
and $(a_{i_2}, b_{j_2}, c_{k_2})\in M_{A,B,C},$ 
we have $a_{i_1}\ne a_{i_2},$ $b_{j_1} \ne b_{j_2},$ and $c_{k_1}\ne c_{k_2}.$ 

\bigskip
\bf
\noindent
Definition 2.5.
\rm\quad
A matching $M_{A,B,C}$ is called perfect if the number of triples in $M_{A,B,C}$ is $n=|A|=|B|=|C|.$

\bigskip
We assume a nonnegative weight of each of the edges is defined for each edge
$w(a_i,b_j),$ $w(b_i,c_j),$ 
and $w(c_i,a_j)$ for any $1\le i\le n$ and $1\le j\le n.$

\bigskip
\bf
\noindent
Definition 2.6.
\rm\quad
If a perfect matching $M_{A,B,C}$ is given, we define its weight $w(M_{A,B,C})$ as
$$w(M_{A,B,C})=\sum_{(a,b,c)\in M_{A,B,C}}\big( w(a,b)+w(b,c)+w(c,a)\big).$$

\bigskip
This definition is in accordance with Definition 1.1 where all the weights
of edges in a complete graph with vertices $a, b,$ and $c$ are included in the sum.

\bigskip
\bf
\noindent
Definition 2.7.
\rm\quad
A perfect matching $M_{A,B,C}^{min}$ is called minimal if its weight is minimal,
that is, $w(M_{A,B,C}^{min})\le w(M_{A,B,C})$ for any other perfect matching $M_{A,B,C}.$

\bigskip
\bf
\noindent
Theorem 2.3.
\rm

Let $A,$ $B,$ and $C$ be the vertex sets of the same cardinality of a complete
tripartite graph $A\cup B\cup C$
with edges in $A\times B,$ $B\times C,$ and $C\times A.$ Then
$$w(M_{A,B}^{min})+w(M_{B,C}^{min})+w(M_{C,A}^{min})\le w(M_{A,B,C}^{min}) .$$

\bigskip
\bf
\noindent
Proof.
\rm\quad
We check that
$$w(M_{A,B}^{min})\le \sum_{(a,b,c)\in M_{A,B,C}^{min}}w(a,b),$$ 
$$w(M_{B,C}^{min})\le \sum_{(a,b,c)\in M_{A,B,C}^{min}}w(b,c),$$ 
$$w(M_{C,A}^{min})\le \sum_{(a,b,c)\in M_{A,B,C}^{min}}w(c,a).$$
Due to definitions 2.2 through 2.7 the sum of these three inequalities 
yields the result.

\bigskip
\bigskip
1) An application of this theorem in a general setting like this
may be found in estimating the accuracy of some heuristic for finding
a perfect matching. 
If we obtain a perfect matching $M_{A,B,C}^{heu}$ in a complete 
tripartite graph, we may use this theorem 2.2 to estimate the accuracy 
of $M_{A,B,C}^{heu}$ as
$$\frac{w(M_{A,B,C}^{heu})}{w(M_{A,B,C}^{min})}\le
\frac{w(M_{A,B,C}^{heu})}{w(M_{A,B}^{min})+w(M_{B,C}^{min})+w(M_{C,A}^{min})}.$$

2) If an inequality in this theorem 2.3 is satisfied as an equality for some perfect 
matching $M_{A,B,C}$, that is, 
$w(M_{A,B,C})=w(M_{A,B}^{min})+w(M_{B,C}^{min})+w(M_{C,A}^{min}),$ we have a minimal solution.

3) The technique of the proof of theorem 2.2 may be used in other situations, such as
$4$-partite matching or $k$-partite matching.

\bigskip
\bf
2.3 Minimal matching on tripartite graphs

\bigskip
\rm
Let $A,$ $B,$ and $C$ be the vertex sets of the same number of vertices. 
We form a complete tripartite graph $A\cup B\cup C$
with edges in $A\times B,$ $B\times C,$ and $C\times A.$

We assume that to each of the vertices in $A,$ $B,$ and $C$ 
real numbers $x_i,$ $y_i,$ and $z_i$ are assigned respectively. 
A weight of each of the edges is defined as
$w(a_i,b_j)=|x_i-y_j|,$ $w(b_i,c_j)=|y_i-z_j|,$ and $w(c_i,a_j)=|z_i-x_j|$ 
for any $1\le i\le n$ and $1\le j\le n.$ Another way to define the weights
is $w(a_i,b_j)=(x_i-y_j)^2,$ $w(b_i,c_j)=(y_i-z_j)^2,$ 
and $w(c_i,a_j)=(z_i-x_j)^2.$ In general the weight has to have property LM.
Without any loss of generality we assume the $n$-tuples $(x_1, x_2,\dots , x_n),$
$(y_1, y_2,\dots , y_n),$ and $(z_1, z_2,\dots , z_n),$ are sorted.
If not, we sort them together with $a_i,$ $b_j,$ and $c_k.$ Obtaining sorted
$n$-tuples can be done in $O(n\log n)$ time.

\bigskip
\bf
\noindent
Theorem 2.4.
\rm

Let $A,$ $B,$ and $C$ be the vertex sets, $|A|=$ $|B|=$ $|C|=n,$ of a complete
tripartite graph $A\cup B\cup C$
with edges in $A\times B,$ $B\times C,$ and $C\times A.$ If the vertices
are assigned real values corresponding to sorted $n$-tuples
$(x_1, x_2,\dots , x_n),$ $(y_1, y_2,\dots , y_n),$ and $(z_1, z_2,\dots , z_n),$
then the minimal matching, with respect to weights with property LM, is given by
$(a_1,b_1,c_1),$ $(a_2,b_2,c_2),$ $\dots ,$
$(a_n,b_n,c_n).$

\bigskip
\bf
\noindent
Proof.
\rm\quad
We claim the matching $(a_1, b_1, c_1),$ $(a_2, b_2, c_2),$
$\dots ,$ $(a_n, b_n, c_n)$ is the minimal one. 
When we form $w(x_i, y_i)+w(y_i, z_i)+w(z_i, x_i)$ 
for each triple separately and add them up over $i$, 
we get the same sum as when we calculate 
$\sum_{i=1}^n w(x_i, y_i)+$ $\sum_{i=1}^n w(y_i, z_i)+$ 
$\sum_{i=1}^n w(z_i, x_i).$ 

It shows we get an equality sign in the inequality in theorem 2.3 which, in turn, means 
that we have obtained a minimal matching.

We would proceed in the same way in the case of weights defined as squares
of differences.

\bigskip
We recall that $|x_i-y_i|+|y_i-z_i|+|z_i-x_i|$ is the distance within this triple
$D(x_i, y_i, z_i)$ introduced in definition 1.1.

\bigskip

\bigskip
\bf
Conclusion

\bigskip
\rm
Results in part 1 may be used as a starting value 
for finding an $n$-matching in a Euclidean  space. 
We may fit a line to data to provide a starting 
$n$-tuple partition followed by a local search. 
One way to do the local search is the concatanation 
of pairs of $n$-tuples to obtain $2n$-tuples and 
enumeration of all the pairs of $n$-tuples. One 
element may be fixed so that we have ${2n-1}\choose{n-1}$ to generate.

The matching algorithm on a line may provide a test
for a general heuristic algorithm for if a general matching 
heuristic works, it should work on a line. A simple heuristic 
may be designed if the vertices are points in a Euclidean 
space, edges are the line segments connecting the vertices, 
and weights are the distances between the end points of 
those line segments. If the number of vertices is $2^n3,$ 
we find the nonbipartite $2$-matching that minimizes the sum 
of the lengths of line segments.
There are $2^{n-1}3$ line segments in this matching. We form
a new graph by taking midpoints of the line segments in the 
matching keeping track of what original vertices the line
segments came from. We repeat this proces until we get three
vertices. Now we work our way back forming a graph with six
vertices and find optimal triples by enumerating all the
pairs of triples of vertices. We continue until we get  
a graph with $2^n3$ vertices. When we use this algorithm
on vertices on a line, we see it gives the correct resullt.

In part 2 of the paper theorem 2.2 may be used in the case
that the weights assigned to edges of a tripartite graph
satisfy the triangle inequality. Let $(a_i, b_j, c_k)$ be given
$a_i\in A,$ $b_j\in B,$ $c_k\in C,$ where $A,B,C$ are disjoint,
$|A|=|B|=|C|=n,$
then $w(c_k,a_i)\le w(a_i,b_j)+w(b_j,c_k).$ 
Let $(a_i, b_j)\in M_{A,B}^{min}$ and $(b_j, c_k)\in M_{B,C}^{min},$
then $(c_k, a_i)$ does not have to be in $M_{C,A}^{min}.$
Matching on a tripartite graph actually asks for 3-cycles 
$a_i,b_j,c_k,a_i.$ 

Actually without knowing or caring what the weights
of $(c_k,a_i)$ are, the use of the triangle inequality gives us 
an upper bound, $w(c_k,a_i)\le$ $w(a_i,b_j)+w(b_j,c_k).$
Thus, if we form a matching like this, denoted as 
$M_{A,B,C}^{\triangle},$ we have 
$$M_{A,B,C}^{\triangle}\le 
2\big(\sum_{(a,b)\in M_{A,B}^{min}}w(a,b)+\sum_{(b,c)\in M_{B,C}^{min}}w(b,c) \big)
=2\big(w(M_{A,B}^{min})+w(M_{B,C}^{min})\big)$$
Thus
$$\frac{w(M_{A,B,C}^{\triangle})}{w(M_{A,B,C}^{min})}\le
\frac{w(M_{A,B,C}^{\triangle})}{w(M_{A,B}^{min})+w(M_{B,C}^{min})+w(M_{C,A}^{min})}\le$$

$$\frac{2\big(w(M_{A,B}^{min})+w(M_{B,C}^{min})\big)}
{w(M_{A,B}^{min})+w(M_{B,C}^{min})+w(M_{C,A}^{min})}\le
\frac{2\big(w(M_{A,B}^{min})+w(M_{B,C}^{min})\big)}
{w(M_{A,B}^{min})+w(M_{B,C}^{min})}=2.$$

\bigskip

\bigskip
\bf
Appendix

\bigskip
\rm
We may try to define a measure of variability in a way different
from the usual approach. We assume there are  $N$ real numbers
$x_1, x_2,\dots , x_N.$ The usual measure of variability, 
the variance $S^2,$ is based on the sum of squares of differences 
from the mean,
$$S^2=\frac{1}{N-1}\sum_{i=1}^N(x_i-\bar{x})^2.$$

 The way we will define the measure of variability without any
 reference to the mean is based on the sum of squares of
 all the differences
 $$\sum_{i=1}^{N-1}\sum_{j=i+1}^N(x_j-x_i)^2.$$
We may check what happens if $y_i=a+bx_i.$

$$\sum_{i=1}^{N-1}\sum_{j=i+1}^N(y_j-y_i)^2=
\sum_{i=1}^{N-1}\sum_{j=i+1}^N(a+bx_j-a-bx_i)^2=
b^2\sum_{i=1}^{N-1}\sum_{j=i+1}^N(x_j-x_i)^2.$$
It means we have the same property for the sum of all the
differences squared and the sum of differences from the mean squared
and it means it is a reasonable characteristic of variability.

Before we show what relation there is between the sum
of squares of all differences and the sum of squares of differences
from the mean we write the sum of squares of all differences as
$$2\sum_{i=1}^{N-1}\sum_{j=i+1}^N(x_j-x_i)^2=
\sum_{i=1}^{N}\sum_{j=1}^N(x_j-x_i)^2.$$

This is easy to see when we write the difference $(x_j-x_i)$
in a different order as $(x_i-x_j).$ When the subscripts 
are the same, we get $x_i-x_i=0.$

Now we review the formula for $(a+b)^2.$ We usually say that
$(a+b)^2=a^2+2ab+b^2$ because we use commutativity
$ab=ba$ therefore $ab+ba=2ab.$ When we don't,
we get $(a+b)^2=aa+ab+ba+bb.$ We will use this idea as

$$(\sum_{j=1}^Nx_j)^2=\sum_{i=1}^N\sum_{j=1}^Nx_ix_j.$$

\bigskip
\bf
\noindent
Theorem
\rm  Let $N>1$ and real numbers $x_1, x_2, \dots , x_N$ be given.
Then
$$\sum_{i=1}^N\sum_{j=1}^N(x_j-x_i)^2=
2N\sum_{i=1}^N(x_i-\bar{x})^2.$$

\bigskip
\bf
\noindent
Proof.
\rm
We expand the formula for twice the sum of squares of all differences

$$\sum_{i=1}^N\sum_{j=1}^N(x_j-x_i)^2=
\sum_{i=1}^N\sum_{j=1}^N(x_j^2+x_i^2-2x_ix_j)=$$
$$\sum_{i=1}^N\sum_{j=1}^Nx_i^2+
\sum_{i=1}^N\sum_{j=1}^Nx_j^2-
2\sum_{i=1}^N\sum_{j=1}^Nx_ix_j=
2N\sum_{i=1}^Nx_i^2-
2\sum_{i=1}^N\sum_{j=1}^Nx_ix_j.$$

We multiply the formula for the sum of squares of differences from the
mean by two

$$2N\sum_{i=1}^N(x_i-\frac{1}{N}\sum_{j=1}^Nx_j)^2=
2N\sum_{i=1}^Nx_i^2-4\sum_{i=1}^Nx_i\sum_{j=1}^Nx_j
+\frac{2}{N}\sum_{i=1}^N(\sum_{j=1}^Nx_j)^2=$$
$$2N\sum_{i=1}^Nx_i^2-4\sum_{i=1}^N\sum_{j=1}^Nx_ix_j
+2\sum_{j=1}^Nx_j)^2=$$
$$2(N\sum_{i=1}^Nx_i^2-4\sum_{i=1}^N\sum_{j=1}^Nx_ix_j
+2\sum_{i=1}^N\sum_{j=1}^Nx_ix_j=$$
$$2N\sum_{i=1}^Nx_i^2-2\sum_{i=1}^N\sum_{j=1}^Nx_ix_j.$$
This proves the desired equality.

\bigskip
\bf
References

\bigskip
\rm
\noindent
Beck, C., Lu, B., Greevy, R. (2015), Nbpmatching: Functions for Optimal Non-Bipartite Matching. R package version 1.4.5. 

\noindent
https://cran.r-project.org/web/packages/nbpmatching

\noindent
http://biostat.mc.vanderbilt.edu/wiki/Main/NonbipartiteMatching

\smallskip
\rm
\noindent
Bondy, J.A., Murty, U.S.R. (1976), Graph Theory with Applications, 
North-Holland, NY.

Garey, M.R., Johnson, D.S. (1979), Computers and Intractability A Guide
to the Theory of NP-Completeness. Freeman and Co, NY.

\smallskip
\rm
\noindent
Garey, M.R., Johnson, D.S. (1979), Computers and Intractability A Guide
to the Theory of NP-Completeness. Freeman and Co, NY.

\smallskip
\rm
\noindent
Lu, B., Greevy, R., Xu, X., Beck, C. (2011),  Optimal Nonbipartite 
Matching and its Statistical Applications. The American Statistician. 
Vol. 65, no. 1, pp. 21-30, .

\smallskip
\noindent
Papadimitriou, C.H., Steiglitz, K. (1982), Combinatorial Optimization:
Algorithms and Complexity, Prentice-Hall, Englewood Cliffs, NJ.

\end{document}